\documentclass[12pt]{article}

 \usepackage {amssymb}
 \usepackage {amsmath}

\usepackage{xcolor}

\def\sqr#1#2{{\vcenter{\vbox{\hrule height.#2pt
              \hbox{\vrule width.#2pt height#1pt \kern#1pt \vrule width.#2pt}
              \hrule height.#2pt}}}}
\def\signed #1{{\unskip\nobreak\hfil\penalty50
              \hskip2em\hbox{}\nobreak\hfil#1
              \parfillskip=0pt \finalhyphendemerits=0 \par}}
\def\endpf{\signed {$\sqr69$}}

\def\dbR{{\mathop{\rm l\negthinspace R}}}

\def\3n{\negthinspace \negthinspace \negthinspace }
\def\2n{\negthinspace \negthinspace }
\def\1n{\negthinspace }

\def\dbE{{\mathop{\rm l\negthinspace E}}}
\def\dbF{{\mathop{\rm l\negthinspace F}}}

\def\ds{\displaystyle}

\def\dbN{{\mathop{\rm l\negthinspace N}}}

\def\dbR{{\mathop{\rm l\negthinspace R}}}
\def\={\buildrel \triangle \over =}

\def\l{\lambda}

\def\th{\theta}

\def\ns{\noalign{\ss} }

\def\no{\noindent}

\def\ms{\medskip}
\def\bs{\bigskip}
\def\q{\quad}
\def\qq{\qquad}

\def\min{\mathop{\rm min}}

\def\pa{\partial}

\def\|{\Big |}
\def\({\Big (}
\def\){\Big )}
\def\[{\Big[}
\def\]{\Big]}
\def\be{\begin{equation}}
\def\bel{\begin{equation}\label}
\def\ee{\end{equation}}
\def\bt{\begin{theorem}}
\def\bcd{\begin{condition}}
\def\ecd{\end{condition}}
\def\et{\end{theorem}}
\def\bc{\begin{corollary}}
\def\ec{\end{corollary}}
\def\bde{\begin{definition}}
\def\ede{\end{definition}}
\def\bl{\begin{lemma}}
\def\el{\end{lemma}}
\def\bp{\begin{proposition}}
\def\ep{\end{proposition}}
\def\br{\begin{remark}}
\def\er{\end{remark}}
\def\ba{\begin{array}}
\def\ea{\end{array}}
\def\ed{\end{document}}
\def\ns{\noalign{\ms}}
\def\ds{\displaystyle}

\def\eps{\epsilon}
\def\square#1{\vbox{\hrule\hbox{\vrule height#1%
     \kern#1\vrule}\hrule}}
\def\rectangle#1#2{\vbox{\hrule\hbox{\vrule height#1%
     \kern#2\vrule}\hrule}}

\font\tenbb=msbm10 \font\sevenbb=msbm7 \font\fivebb=msbm5

\newfam\bbfam
\scriptscriptfont\bbfam=\fivebb \textfont\bbfam=\tenbb
\scriptfont\bbfam=\sevenbb

\newtheorem{lemma}{Lemma}[section]
\newtheorem{remark}{Remark}[section]

\newtheorem{theorem}{Theorem}[section]
\newtheorem{corollary}{Corollary}[section]

\newtheorem{definition}{Definition}[section]
\newtheorem{proposition}{Proposition}[section]
\newtheorem{condition}{Condition}[section]

\makeatletter
   
   \@addtoreset{equation}{section}
\makeatother

\begin{document}
\title{\bf Null controllability for
stochastic semi-discrete parabolic  equations  \thanks{This  work  is  partially  supported  by  the  NSF  of  China  under  grants  11971333 and  11931011. The author gratefully acknowledges Professor Xiaoyu Fu  for her
guidance and suggestions.}}

\author{Qingmei Zhao\thanks{School of Mathematics Sciences, Sichuan Normal University, Chengdu 610066, China. E-mail address: qmmath@163.com.}  }

\date{}

\maketitle

\begin{abstract}
\no  In this paper, we present a null controllability result for a class of  stochastic semi-discrete   parabolic equations. For this purpose, an observability estimate is established for backward stochastic semi-discrete parabolic equations,  with an explicit observability constant that depends on the discretization parameter. This estimate is obtained by a new Carleman estimate for backward stochastic semi-discrete parabolic operators.

 \end{abstract}

\bs

\no{\bf 2010 Mathematics Subject Classification}.  Primary 93B05; Secondary 93B07, 93C20
\bs

\no{\bf Key Words}.  Stochastic semi-discrete  parabolic equations,  controllability, observability,  global Carleman
estimate.

\ms

\section{Introduction and main result}
Let $T>0$ and $(\Omega, \mathcal{F}, \{\mathcal{F}_t\}_{t\ge
0}, \mathcal{P})$ be a complete filtered probability space on which a one-dimensional
standard Brownian motion  $\{B(t)\}_{t\geq 0}$ is defined such that
$\{\mathcal{F}_t\}_{t\geq 0}$ is the natural filtration generated by
$B(\cdot)$, augmented by all the $\mathcal{P}$-null sets in
$\mathcal{F}$. Write $\mathbb{F}$ for the progressive $\sigma$-field with respect to $\mathcal{F}$. Let $\mathcal{H}$ be a Banach space. Denote by $L^2_{\mathcal{F}_\tau}(\Omega; \mathcal{H})$ the space of all $\mathcal{F}_\tau$-measurable random variables $\xi$ such that $\dbE | \xi |^2_{\mathcal{H}}<\infty$; by $L^2_{\mathbb{F}}(0,T;\mathcal{H})$ the
Banach space consisting of all $\mathcal{H}$-valued
$\{\mathcal{F}_t\}_{t\geq 0}$-adapted processes $X(\cdot)$ such that
$\mathbb{E}(|X(\cdot)|^2_{L^2(0,T;\mathcal{H})})<\infty$, with the
canonical norm; by $L^\infty_{\mathbb{F}}(0,T;\mathcal{H})$ the
Banach space consisting of all $\mathcal{H}$-valued
$\{\mathcal{F}_t\}_{t\geq 0}$-adapted essentially bounded processes;
and by $L^2_{\mathbb{F}}(\Omega; C([0,T];\mathcal{H}))$ the Banach
space consisting of all $\mathcal{H}$-valued
$\{\mathcal{F}_t\}_{t\geq 0}$-adapted continuous processes
$X(\cdot)$ such that
$\mathbb{E}(|X(\cdot)|^2_{C([0,T];\mathcal{H})})<\infty$.

Let $\omega$  be a given
nonempty open subset of $(0, 1)$, and denote by $\chi_{\omega}$ the
characteristic function of the set $\omega$. We recall the following controlled stochastic parabolic equation  within a continuous framework:
\begin{eqnarray}\label{e*01}\left\{
\begin{array}{lll}
\ds dy- y_{xx}dt
=(a_1 y+\chi_\omega u)dt+(a_2y+v)dB(t)  &\mbox{ in }(0, 1)\times(0,T),\\[3mm]
y(0,t)=0, \ y(1,t)=0& \mbox{ on } (0,T),\\[3mm]
 y(x, 0)=y_0(x) &\mbox{ in }(0, 1),
\end{array}
\right.
\end{eqnarray}
with suitable coefficients $a_1(\cdot), a_2(\cdot)\in  L^\infty_{\mathbb{F}}(0, T; L^\infty(0,1))$. In (\ref{e*01}), the control variable consists of the pair $(u, v)$, $y$ is the state, and $y_0\in L^2(\Omega, \mathcal{F}_0, \mathcal{P}; L^2(0,1))$ is the initial datum. The system
(\ref{e*01}) is said to be null controllable at time $T$ if for any given $y_0$, there exist controls $(u,v)\in  L^2_{\mathbb{F}}(0, T; L^2(\omega))\times L^2_{\mathbb{F}}(0, T; L^2(0,1))$ such
that  the corresponding solution satisfies $y(x, T) = 0$, $\mathcal{P}\mbox-{a.s.}$

In the literature, significant progress has been achieved in investigating null controllability for both deterministic and stochastic parabolic equations within the continuous framework(see \cite{FI, FL, FLZ, LZ1, TZ} and the references therein). More recently,  important progress  has been made in addressing control problems for semi-discrete approximations of deterministic parabolic equations (see \cite{AB, CGL, FM, FM1, MZ} and the references therein). However, to the best of our knowledge, within the stochastic framework, there is a noticeable gap in addressing the corresponding controllability issues for stochastic semi-discrete parabolic equations.

For given $N \in \dbN$, we
set the space discretization parameter $h := 1/(N + 1)$. We consider the pairs $(x_i,t)$ with $t\in(0, T)$,
and $x_i = ih, i = 0, ...,N+1$. Applying the centered finite difference method to the spatial variable for the
system (\ref{e*01}), %and considering two additional points $x_{-1} := x_0- h$ and $x_{N+2} := x_{N+1} + h$ to discretize the boundary conditions,
we obtain the following semi-discrete system:
\begin{eqnarray}\label{e*03}\left\{
\begin{array}{lll}
\ds dy^i- \frac{1}{h^2}(y^{i+1}-2y^i+y^{i-1})dt
=(a^i_1 y^i+\chi_\omega^i u^i)dt+(a^i_2y^i+v^i)dB(t), \\[3mm]
\qq\qq\qq\qq\qq\qq\qq \q t\in (0, T),\
 \ \ds i = 1, ..., N,\\[3mm]
y^0(t)=0, \ y^{N+1}(t)=0,\ \  t\in (0,T),\\[3mm]
 y^i(0)=y^i_0, \ \ i =  1, ..., N,
\end{array}
\right.
\end{eqnarray}
where
$a^i_1=a_1(x_i)$, $a^i_2=a_2(x_i)$, $\chi_\omega^i =\chi_\omega(x_i)$, $u^i =u(x_i)$, and $v^i =v(x_i)$.

To present our main result, we first introduce some notations. First, we  define a regular partition of the interval  $[0, 1]$ as follows:
$$
\mathcal{K}:=\{x_i\ | \ i =0,1, ..., N,N+1\}.
$$
Considering any set of points $\mathcal{W}\subset \mathcal{K}$, we define the dual meshes $\mathcal{W}'$ and
$\mathcal{W}^*$, respectively, by
\begin{eqnarray}\label{mesh1}
\mathcal{W}':=\tau_+(\mathcal{W})\cap\tau_-(\mathcal{W}),\ \ \mathcal{W}^*:=\tau_+(\mathcal{W})\cup\tau_-(\mathcal{W}),
\end{eqnarray}
where
$$
\tau_\pm(\mathcal{W}):=\{x\pm\frac{h}{2}\ | \ x\in\mathcal{W}\}.
$$
We denote $\overline{\mathcal{W}}:=(\mathcal{W}^*)^*$ and $\mathring{\mathcal{W}}:=(\mathcal{W}')'$.
Note that if $\mathring{\overline{\mathcal{W}}} = \mathcal{W}$, then for two consecutive points
$x_i, x_{i+1} \in \mathcal{W}$, we have $x_{i+1} - x_i = h$. Thus, any subset $\mathcal{W}\subset \mathcal{K}$ that verifies $\mathring{\overline{\mathcal{W}}} = \mathcal{W}$ is called regular mesh.
Finally, we define the boundary of a regular mesh $ \mathcal{W}$ as $\partial{ \mathcal{W}} := \overline{ \mathcal{W}}\backslash\mathcal{W}$.

We introduce, using (\ref{mesh1}), the semi-discrete sets. We define the semi-discrete set
$G := \mathcal{W} \times (0, T)$. Moreover, we say that the semi-discrete set is regular if the space variable is a regular
mesh. Additionally, we define the dual semi-discrete sets by $G' := \mathcal{W}' \times(0, T)$ and $G^* := \mathcal{W}^*\times(0, T)$. Similarly, the
semi-discrete boundary is given by $\partial G = \partial{ \mathcal{W}} \times (0, T)$.

Next, we define the average operator $A_h$ and the difference operator $D_h$  on $G$ by
\begin{eqnarray*}
A_h(u)(x,t):=\frac{\tau_+u(x,t)+\tau_-u(x,t)}{2},\\
D_h(u)(x,t):=\frac{\tau_+u(x,t)-\tau_-u(x,t)}{h},
\end{eqnarray*}
where
$\tau_\pm u(x,t)=u(x\pm\frac{h}{2},t)$.

Here and in what follows, denote by $C(\mathcal{W})$  the set of real-valued functions defined in $\mathcal{W}$,  by $L^2_h(\mathcal{W})$ the Hilbert space with
 inner product
 $$
 \langle u,v\rangle_\mathcal{W}:=\int_\mathcal{W}u v:=h\sum_{x\in\mathcal{W}}u(x)v(x).
 $$
Further, we define the outward normal for $(x,t) \in\partial G$ as
\begin{eqnarray*}\nu(x,t):=\left\{
\begin{array}{llr}
\ds 1,  &\ if \  (\tau_-(x),t)\in G^*\  and \ (\tau_+(x),t)\notin G^*,\\[3mm]
-1,   &\ if \  (\tau_-(x),t)\notin G^*\  and \ (\tau_+(x),t)\in G^*,\\[3mm]
 0, &\  otherwise.
\end{array}
\right.
\end{eqnarray*}

Meanwhile, we introduce the trace operator on $\partial G$ as
\begin{eqnarray*}\forall (x,t)\in \partial G,\ t_r(u):=\left\{
\begin{array}{llr}
\ds \tau_-u(x,t),  &\ if \  \nu(x,t)=1,\\[3mm]
\tau_+u(x,t),  &\ if \  \nu(x,t)=-1,\\[3mm]
 0, & \ if \   \nu(x,t)=0.
\end{array}
\right.
\end{eqnarray*}
Finally, we  introduce the discrete integration on the boundary on $\partial \mathcal{W}$ as 
$$
\int_{\partial  \mathcal{W}}u:=\sum_{x\in\partial\mathcal{W}}u(x).
$$
Thus, defining $Q := \mathcal{M}\times(0, T)$ where $\mathcal{M}:= \mathring{\mathcal{K}}$, the controlled semi-discrete
system (\ref{e*03}) can be written as
\begin{eqnarray}\label{e*03a}\left\{
\begin{array}{lll}
\ds dy- D_h^2ydt
=(a_1 y+\chi_\omega v)dt+(a_2y+v)dB(t),
\  (x,t)\in Q,\\[3mm]
y(0,t)=0, \ y(1,t)=0,\ \  t\in (0,T),\\[3mm]
 y(x,0)=y_0, \ x\in \mathcal{M}.
\end{array}
\right.
\end{eqnarray}
The main purpose of this paper is to study the null controllability for the semi-discrete stochastic parabolic equation (\ref{e*03a}).
Similar to the continuous case in (\ref{e*01}), for any given initial value $y_0\in L^2(\Omega,\mathcal{F}_0, \mathcal{P};  L^2_h(\mathcal{M}))$, we wonder to find controls $(u, v) \in L^2_{\mathbb{F}}(0, T; L^2_h(\omega\cap\mathcal{M}))\times L^2_{\mathbb{F}}(0, T; L^2_h(\mathcal{M}))$
such the solution to system (\ref{e*03a}) satisfies $y_i(T) = 0$, $\mathcal{P}\mbox-{a.s.}$,  for all $i = 1, ..., N$. If the control exists
we would say that the system (\ref{e*03a}) is uniformly null-controllable.

It is well known that the  uniform null controllability
property is very hard to address even for deterministic case. The underlying difficulty lies in demonstrating  some uniform behavior with respect to the discretization parameter.
For the uniform null controllability of deterministic semi-discrete parabolic equations,
there are some results in some special cases (see \cite{LT, LZ, Zua1, Zua2}).
However, as pointed out in \cite{ Zua1}, these methods are not applicable  to  semi-discrete parabolic equations with variable coefficients.
Uniform null controllability
property may not even hold in the two-dimensional case, as demonstrated by the counterexample provided in \cite{ Zua2}. Consequently, alternative definitions have been proposed to relax the requirement for uniform null controllability. One such relaxed definition is  $\phi$-null controllability, i.e. $|y(T)|_{L_h^2(\mathcal{M})}\le
\phi(h)$, where $\phi$ is a real-valued function that tend to zero when space
discretization parameter $h$ tends to zero. This definition proves to be highly useful for obtaining well-characterized controls.
%We found the papers  [] where different methods
%are proposed to deal with this $\phi$-null controllability problem for semi-discrete parabolic equations in deterministic case.

On the other hand, in practice, due to the interference of random factors, stochastic processes give a natural replacement for deterministic functions in mathematical descriptions. Compared with deterministic case, some substantially difficulties  arise in the study of the stochastic partial differential equations.
For example, the solution to a stochastic partial differential equation is non-differentiable with respect to noise variable, and the usual compactness embedding result is not valid for solution spaces of the stochastic evolution equations.
Indeed, many tools and methods, which are effective in the deterministic case, do not work
anymore in the stochastic setting.

Despite these challenges, to the best of our knowledge, no research has been conducted on the null controllability of stochastic semi-discrete parabolic equations. In this paper, we address this gap by investigating the controllability problems of general stochastic semi-discrete parabolic equations.  To achieve this objective, we establish an observability estimate for backward stochastic semi-discrete parabolic equations, with an explicit observability constant that depends on the discretization parameter. This estimate is derived through a novel Carleman estimate tailored for backward stochastic semi-discrete parabolic operators. It is well known that Carleman estimate is a powerful tool in the study of control problems and inverse problems for partial differential equations. Here, we refer to \cite{FI} for the Carleman estimate of deterministic parabolic operator in continuous framework, \cite{TZ} for the Carleman estimate of continuous stochastic parabolic operator, \cite{BHR, BHR1, BR, CLNP, N} for Carleman estimates  of deterministic semi-discrete elliptic equations and parabolic equations. Unlike the references mentioned above, establishing Carleman estimates for stochastic semi-discrete parabolic operators poses new challenges. To address this, we introduce a new auxiliary function to tackle the difficulties arising from the diffusion term.

Our main result of this paper can be stated as follows.

\bt
There exist $C$ and $h_0$, for all $h\le h_0$, there exist $(u, v) \in L^2_{\mathbb{F}}(0, T; L^2_h(\omega\cap\mathcal{M}))\times L^2_{\mathbb{F}}(0, T; L^2_h(\mathcal{M}))$ such that the solution to $(\ref{e*03a})$ satisfies
\begin{eqnarray*}
\begin{array}{rl}
&\ds \dbE\int_{Q}|v|^2dt+\dbE\int_{0}^T\int_{\omega\cap\mathcal{M}}|u|^2 dt
 \le  C\dbE\int_{\mathcal{M}}|y_0|^2,
 \end{array}
\end{eqnarray*}
and
\begin{eqnarray*}
\begin{array}{rl}
&\ds \dbE \int_{\mathcal{M}}|y(T)|^2
 \le   Ce^{-\frac{C}{h}}\dbE\int_{\mathcal{M}}|y_0|^2.
 \end{array}
\end{eqnarray*}
\et

The rest of this paper is organized as follows. In section 2, we give some needed preliminaries.
In section 3, we prove a global Carleman estimate for a backward  stochastic semi-discrete parabolic operator. In section 4, we prove the
observability result for the backward stochastic semi-discrete parabolic equation. Finally, in section 5, we prove the main controllability result.

\section{Preliminaries}

In this section, we give some preliminaries we needed.

\begin{lemma}\label{l1}
\rm{(\cite{CLNP})} For any $u,v\in C(\overline{\mathcal{W}})$, we have for the difference operator
\begin{eqnarray}\label{p1}
D_h(uv)=D_huA_hv+A_huD_hv,\  on\  \mathcal{W}^*.
\end{eqnarray}
Similarly, the average of the product gives
\begin{eqnarray}\label{p2}
A_h(uv)=A_huA_hv+\frac{h^2}{4}D_huD_hv,\  on\  \mathcal{W}^*.
\end{eqnarray}
Finally, on $\mathring{\mathcal{W}}$ we have
\begin{eqnarray}\label{p3}
u=A^2_hu-\frac{h^2}{4}D^2_hu.
\end{eqnarray}
\end{lemma}
Next, we recall the following integration by parts formulas for the discrete difference and average operators.
\begin{lemma}\label{l2}
\rm{(\cite{CLNP})} Let $\mathcal{W}$ be a semi-discrete regular mesh. For $u\in C(\overline{\mathcal{W}})$ and $v\in C(\mathcal{W}^*)$, we have
\begin{eqnarray}\label{p4}
\int_\mathcal{W}uD_hv=-\int_{\mathcal{W}^*}D_hu v+\int_{\partial \mathcal{W}}u t_r(v)\nu
\end{eqnarray}
and
\begin{eqnarray}\label{p5}
\int_\mathcal{W}uA_hv=\int_{\mathcal{W}^*}A_hu v-\frac{h}{2}\int_{\partial \mathcal{W}}u t_r(v).
\end{eqnarray}
\end{lemma}

Further, we have the following Lemma, which gives the relationship  between disctete operators and derivative operators.

\begin{lemma}\label{l5}
\rm{(\cite{CLNP})} Let $f$ be a $(n+4)$-times differentiable function defined on $\dbR$ and $m,n\in\dbN$. Then
\begin{eqnarray*}
\begin{array}{rl}
&A_h^mD_h^n f=f^{(n)}+R_{D_h^n}(f)+R_{A_h^m}(f)+R_{A_h^mD_h^n}(f),
\end{array}
\end{eqnarray*}
where
$$
R_{D_h^n}(f):=h^2\sum_{k=0}^n \binom{n}{k}(-1)^k\(\frac{n-2k}{2}\)^{n+2}\int_0^1\frac{(1-\sigma)^{n+1}}{(n+1)!} f^{(n+2)}\(\cdot+\frac{(n-2k)h}{2}\sigma\)d\sigma,
$$
$$
R_{A_h^m}(f):=\frac{h^2}{2^{m+2}}\sum_{k=0}^n \binom{m}{k}(m-2k)^2\int_0^1(1-\sigma) f^{(2)}\(\cdot+\frac{(m-2k)h}{2}\sigma\)d\sigma,
$$
and $R_{A_h^mD_h^n}=R_{A_h^m}\circ R_{D_h^n}$.
\end{lemma}

In what follows, we will present some results related to discrete operations performed on the Carleman weight functions. To begin with, we
introduce the weight functions. By \cite{FI}, we know  that there exists a real-valued
function $\psi\in C^p(\widetilde G)$, $p$ sufficiently large,  so that
\bel{240107-0a} \psi(x)> 0 \ \text{in} \ \widetilde G; \ \   |\partial_x\psi(x) | > 0 \ \text{in} \  \widetilde G \backslash \omega_0; \ \  \text{and} \ (\partial_x\psi)(0)>0, \ (\partial_x\psi)(1)<0,\ee
where $\widetilde G$ is a smooth open and connected neighborhood of $[0,1]$ and $\omega_0$ is any given nonempty open subset of $\omega$ satisfying $ \overline{\omega_0}\subseteq \omega$.

For any positive parameters $\l>1$, $\mu>1$, write
\begin{eqnarray}\label{***f}
\begin{array}{rl}
&\phi(x)=e^{\mu \psi(x) }-e^{2\mu |\psi|_{\infty}}, \ \varphi(x)=e^{\mu \psi(x) },\   r(t,x)=e^{s(t)\phi(x)},\ \rho(t, x)=r^{-1},\\[5mm]
&\ds  s(t)=\l \th(t), \ \theta(t)=\displaystyle{\frac{1}{(t+\delta T)(T+\delta T-t)}},
\end{array}
\end{eqnarray}
for $0< \delta<\frac{1}{2}$.

In the sequel, for any $n\in\dbN$, we denote by $\mathcal{O}(s^n)$ a function of
order $s^n$, for sufficiently large $s$, and by $\mathcal{O}_\mu(s^n)$ a function of order $n$
for fixed $\mu$ and sufficiently
large $s$.

By (\ref{***f}) and Lemma \ref{l5}, we have the following estimation of coefficients involved in the energy terms of Carleman estimate.

\begin{lemma}\label{l3} For $\l h (\delta T^2)^{-1}\le 1$, we have
\bel{240107-a}\left\{
\ba{lll}
&rD_h^2\rho=A_h(r D_h^2\rho)=\mathcal{O}_\mu(s^2) ,&\partial_t(r D_h^2\rho)=\partial_tA_h(r D_h^2\rho)=s^2 T\th \mathcal{O}_\mu(1),\\\ns
&\ds D^2_h(rD_h^2\rho)=\mathcal{O}_\mu(s^2) ,&\partial_tD^2_h(rD_h^2\rho)=s^2 T\th \mathcal{O}_\mu(1),\\\ns
&\ds A_hD_h(rD_h^2\rho)=\mathcal{O}_\mu(s^2),
\ea\right.
\ee
and
\bel{240107-b}\left\{
\ba{ll}
&r^2 D_h^2\rho A_hD_h\rho=\mathcal{O}_\mu(s^3),\q D_h(r^2 D_h^2\rho A_hD_h\rho)=\mathcal{O}_\mu(s^3),\\
\ns & A_hD_h(r^2 D_h^2\rho A_hD_h\rho)
=-3s^3\varphi^3\mu^4(\partial_x\psi)^4+s^3\varphi^3\mathcal{O}(\mu^3)+\mathcal{O}_\mu(s^2)+s^3\mathcal{O}((sh)^2).
\ea\right.
\ee
And moreover,
\bel{240107-c}\left\{
\ba{ll}
&r^2 A_h^2\rho A_hD_h\rho=-s\varphi\mu \partial_x \psi +s\mathcal{O}_\mu((sh)^2),\\
\ns &D_h(r^2 A_h^2\rho A_hD_h\rho)=-s\varphi \mu^2(\partial_{x}\psi)^2+s\varphi\mathcal{O}(\mu)+s\mathcal{O}_\mu((sh)^2).
\ea\right.
\ee
\end{lemma}

\noindent{\bf  Proof. } By (\ref{***f}) and Lemma \ref{l5}, one can obtain the  estimations of (\ref{240107-a})--(\ref{240107-c}) directly. For the reader's convenience, here we only give the estimations of $rD_h^2\rho$, $\partial_tA_h(r D_h^2\rho)$, $A_hD_h(r^2 D_h^2\rho A_hD_h\rho)$ and $D_h(r^2 A_h^2\rho A_hD_h\rho)$.

By (\ref{***f}), it is easy to see that
 \bel{240107-d}\left\{\ba{ll}
 &\ds |\partial _t\th|\le CT\th^2, \q |\partial _{tt}\th|\le CT\th^3,\\
  \ns&\ds \pa_x \rho=-s \mu \varphi\rho\partial_{x}\psi, \q\ \pa_x^2\rho=s^2\mu^2\varphi^2\rho(\partial_{x}\psi)^2+s\varphi\rho\mathcal{O}(\mu^2).
 \ea\right.\ee
Then, by Lemma \ref{l5}, we have
\begin{eqnarray*}
\begin{array}{rl}
&rD_h^2\rho=r\partial_x ^2 \rho+rR_{D_h^2}(\rho)=s^2\varphi^2\mathcal{O}(\mu^2)+s^2\mathcal{O}_\mu((sh)^2)=\mathcal{O}_\mu(s^2),
\end{array}
\end{eqnarray*}
which is used the fact that $sh\le 1$ and $r(x) \partial_x ^4\rho(x+\sigma h)=\mathcal{O}_\mu(s^4)$ for any $\sigma\in[0, 1]$.
By Lemma \ref{l5}, it follows that
$$
A_h(r D_h^2\rho)=A_h(r\partial_x ^2 \rho)+A_h[rR_{D_h^2}(\rho)]=r\partial_x ^2 \rho+R_{A_h}(r\partial_x ^2 \rho)+rR_{D_h^2}(\rho)+R_{A_h}[rR_{D_h^2}(\rho)].
$$
For $\sigma\in[0, 1]$, $\partial_t\partial_x ^2(r\partial_x ^2 \rho)(x+\sigma h, t)=T\th  \mathcal{O}_\mu(s^2)$ by (\ref{240107-d}), then $\partial_t[ R_{A_h}(r\partial_x ^2 \rho)]=T\th\mathcal{O}_\mu((sh)^2)$. And by $\partial_t[r(x, t) \partial_x ^4\rho(x+\sigma h, t)]=T\th  \mathcal{O}_\mu(s^4)$, we have $\partial_t[ rR_{D_h^2}(\rho)]=s^2 T\th \mathcal{O}_\mu((sh)^2)$. Similarly, $\partial_t[R_{A_h}[rR_{D_h^2}(\rho)]=T\th  \mathcal{O}_\mu((sh)^4)$. Then, for $sh\le 1$, it holds that
 \begin{eqnarray*}
\begin{array}{rl}
&\partial_tA_h(r D_h^2\rho)=T\th  \mathcal{O}_\mu(s^2)+s^2 T\th \mathcal{O}_\mu((sh)^2)+T\th  \mathcal{O}_\mu((sh)^4)=T\th  \mathcal{O}_\mu(s^2).
\end{array}
\end{eqnarray*}
Next, by Lemma \ref{l5}, we have
$$
 D_h^2\rho = \partial_x^2\rho +R_{D_h^2}(\rho),\q  A_hD_h\rho=\partial_x \rho+R_{A_h}(\rho)+R_{D_h}(\rho)+R_{D_hA_h}(\rho).
$$
Thus, combining the the above estimate, Lemma \ref{l5} and  (\ref{240107-d}), we conclude that for $sh\le1$,
\begin{eqnarray*}
\begin{array}{rl}
 A_hD_h(r^2 D_h^2\rho A_hD_h\rho)&=\partial_x(r^2 \partial_x^2\rho \partial_x \rho)+s^3\mathcal{O}_\mu((sh)^2)\\\ns
&=\partial_x[-s^3\mu^3\varphi^3(\partial_{x}\psi)^3]+s^2\varphi^2\mathcal{O}(\mu^4)+s^3\mathcal{O}((sh)^2)\\\ns
& =-3s^3\varphi^3\mu^4(\partial_x\psi)^4+s^3\varphi^3\mathcal{O}(\mu^3)+s^2\varphi^2\mathcal{O}(\mu^4)+s^3\mathcal{O}((sh)^2).
\end{array}
\end{eqnarray*}
 Similarly, by Lemma \ref{l5} and (\ref{240107-d}), for $sh\le1$, we can get that
\begin{eqnarray*}
\begin{array}{rl}
D_h(r^2 A_h^2\rho A_hD_h\rho)&=\partial_x(r^2 \rho \partial_x \rho)+s\mathcal{O}_\mu((sh)^2)\\\ns
&=\partial_x(-s \mu \varphi\partial_{x}\psi)+s\mathcal{O}_\mu((sh)^2)\\\ns
&\ds=-s\varphi \mu^2(\partial_{x}\psi)^2+s\varphi\mathcal{O}(\mu)+s\mathcal{O}_\mu((sh)^2).
\end{array}
\end{eqnarray*}
This completes the proof of Lemma \ref{l3}. \endpf

\ms
\br
In fact, what is needed in the proof of Lemma  $\ref{l3}$ is $sh\le 1$. Notice that $s=\l\th$ and $\th(t)\le\th(0)=\th(T)< (\delta T^2)^{-1}$, for all $t\in [0, T]$. Therefore, the condition we present in Lemma $\ref{l3}$ is $\l h (\delta T^2)^{-1}\le 1$.
\er

Finally, for any $v\in L^2_{\mathbb{F}}(\Omega; C([0,T]; L^2_h(\mathcal{W})))$,  we use the following notations:
\begin{eqnarray}\label{x1}
\begin{array}{rl}
& X^1:=\ds\dbE\int_{{Q}}\Big[s^3\varphi^3\mathcal{O}(\mu^3)+ \mathcal{O}_\mu(s^2)+ \th\mathcal{O}_\mu(s^2)+s^3\mathcal{O}_\mu((sh)^2)\Big]|v|^2dt\\\ns
&\ds\q \q\q
+\dbE\int_{Q^*}\Big[s\varphi\mathcal{O}(\mu)+ \mathcal{O}_\mu(1)+\th\mathcal{O}_\mu(1)+s\mathcal{O}_\mu((sh)^2)) \Big]|D_hv|^2dt,
\end{array}
\end{eqnarray}
\begin{eqnarray}\label{x2}
\begin{array}{rl}
& X^2:=\ds\dbE\int_{{Q}}\mathcal{O}_\mu(s^2)|dv|^2
+\dbE\int_{Q^*} \mathcal{O}_\mu((sh)^2)|D_h(dv)|^2,
\end{array}
\end{eqnarray}
\begin{eqnarray}\label{x3}
\begin{array}{rl}
& X^3:=\ds
\dbE\int_{\partial Q}\Big[s\mathcal{O}_\mu((sh)^2)+\mathcal{O}_\mu(1)\Big]t_r(|D_hv|^2)dt,
\end{array}
\end{eqnarray}
and
\begin{eqnarray}\label{x4}
\begin{array}{rl}
& X^4:=\ds\dbE\int_{{ \mathcal{M}}}\mathcal{O}_\mu(s^2)v^2\Big|_{t=0}
 \ds+\dbE\int_{\mathcal{M}^*} \Big[\mathcal{O}_\mu(1)+\mathcal{O}_\mu((sh)^2)\Big]|D_hv|^2\Big|_{t=0}\\\ns
&\q\q\q \ds+\dbE\int_{{ \mathcal{M}}}\mathcal{O}_\mu(s^2)v^2\Big|_{t=T}
 \ds+\dbE\int_{\mathcal{M}^*} \Big[\mathcal{O}_\mu(1)+\mathcal{O}_\mu((sh)^2)\Big]|D_hv|^2\Big|_{t=T}.
\end{array}
\end{eqnarray}
In fact, we can see that $X^1$ is the lower order term about energy, $X^2$ is the estimate about  drift term, $X^3$ is the lower order term about space boundary, and $X^4$ is the time boundary term.

\section{Semi-discrete Carleman estimate for uniform meshes}
In the section, we give the semi-discrete Carleman estimate for the backward stochastic  semi-discrete  parabolic operator.
 Then, we have the following Carleman estimate.

\bt
There exists a positive constant $\mu_0$ such that for any $\mu>\mu_0$, one can find two positive
constants $\varepsilon_0>0$ and $h_0>0$, with $0<\varepsilon_0\le 1$ and a sufficiently large $\l_0\ge 1$, such that for $\l\ge \l_0$, and $h\le h_0 $ and $\l h(\delta T^2)^{-1}\le \varepsilon_0$, it holds that
\begin{eqnarray*}
\begin{array}{rl}
&\ds\dbE\int_{Q}  s^3e^{2s\phi} w^2dt+\ds\dbE\int_{Q^*}se^{2s\phi} |D_hw|^2dt\\\ns
&\le\ds C\(\dbE\int_0^T\int_{\omega\cap \mathcal{M}}s^3e^{2s\phi} w^2dt
+ \dbE\int_{{Q}}s^2e^{2s\phi}|g|^2
 +\dbE\int_{Q}e^{2s\phi} | f|^2dt\\\ns
&\ds\q +h^{-2}\dbE\int_{\mathcal{M}}e^{2s\phi} |w|^2\Big|_{t=0}
+h^{-2}\dbE\int_{\mathcal{M}} e^{2s\phi}  |w|^2\Big|_{t=T}\),
\end{array}
\end{eqnarray*}
for all $w\in L^2_{\mathbb{F}}(\Omega; C([0,T]; L^2_h(\mathcal{W}))) $ satisfying $dw+ D_h^2wdt=fdt+gdB(t)$ with $w=0$ on $\partial {\mathcal{M}}$.
\et

\noindent{\bf  Proof. }  The whole proof is divided  into six steps.

\ms

\ms

\noindent{\bf Step 1.}
At first, let
\begin{eqnarray}\label{f}
\begin{array}{rl}
\mathcal{L}_h w:=\mbox{d}w+ D_h^2wdt=fdt+gdB(t).
\end{array}
\end{eqnarray}
Setting $v = rw$  and noting that $\rho=r^{-1}$, we have
\begin{eqnarray}\label{car01}
r\mathcal{L}_h(w)=r\mathcal{L}_h(\rho v)=I_1dt-\Phi(v)dt+I_2,
\end{eqnarray}
where
\begin{eqnarray}\label{car02}
\left\{\begin{array}{rl}
&\ds I_1=r D_h^2\rho A_h^2v+r A_h^2\rho D_h^2v-\lambda\phi\partial _t\th v+\Phi(v)\ds,\\ \ns
&I_2=dv+2r A_hD_h\rho A_hD_hvdt,\\\ns
&\ds\Phi(v)=\frac{h^2}{2}D_h\[D_h(r D_h^2\rho)A_hv\].
\end{array} \right.\end{eqnarray}
Therefore, by (\ref{car01}), it holds that
\begin{eqnarray}\label{car03}
2r\mathcal{L}_h(w)I_1=2|I_1|^2dt-2I_1\Phi(v)dt+2I_1I_2.
\end{eqnarray}
Hence, combining (\ref{car03}) with (\ref{f}), it follows that
\begin{eqnarray*}
&\ds\dbE\int_Q2r  fI_1dt=2\dbE\int_Q|I_1|^2dt-2\dbE\int_QI_1\Phi(v)dt+ 2\dbE\int_Q I_1I_2.
\end{eqnarray*}
By using the triangle inequality, it is easy to see that
\begin{eqnarray}\label{car04}
\dbE\int_0^T|r f|^2dt\ge 2\dbE\int_Q I_1I_2-\dbE\int_Q|\Phi(v)|^2dt.
\end{eqnarray}
The next step is to provide an estimate for the right-hand side of (\ref{car04}).
 Denote the terms in the  right-hand side of $I_1$ and $I_2$ in (\ref{car02})
by $I_1^i (i = 1, 2, 3, 4)$ and $I_2^j (j = 1, 2)$, respectively. Then
\begin{eqnarray}\label{car05}
 2\dbE\int_Q I_1I_2= \sum_{i=1}^4\sum_{j=1}^2 J_{ij},\q J_{ij}:=2\dbE\int_QI_1^iI_2^j.
\end{eqnarray}

\medskip

\noindent{\bf Step 2. } Let us compute ``$J_{11}+J_{21}+J_{31}+J_{41}$".
 By (\ref{p3}) and (\ref{car02}), noting that $r\rho=1$, we have
\begin{eqnarray}\label{car011}
\begin{array}{rl}
&J_{11}+J_{21}+J_{31}+J_{41}\\\ns
&\ds=2\ds\dbE\int_Q\Big(r D_h^2\rho A_h^2v+ r A_h^2\rho D_h^2v-\lambda\phi\partial _t\th v+\Phi(v)\Big)dv
\\\ns
&\ds  =2\dbE\int_{Q}\Big(r D_h^2\rho v+\frac{h^2}{4}r D_h^2\rho D_h^2v+ r \rho D_h^2v+\frac{h^2}{4}r D_h^2\rho D_h^2v-\lambda\phi\partial _t\th v+\Phi(v)\Big)dv\\\ns
&\ds  =2\dbE\int_{Q}\Big(q^{11} v-\lambda\phi\partial _t\th v+ D_h^2v+\frac{h^2}{2}q^{11} D_h^2v+\Phi(v)\Big)dv,
\end{array}
\end{eqnarray}
where
 \bel{24-05-0a}
  q^{11}:=r D_h^2\rho.
  \ee

For $\l h (\delta T^2)^{-1}\le 1$, by (\ref{240107-a}), (\ref{240107-d}), (\ref{24-05-0a}) and It\^o's formula, we can get that
 \begin{eqnarray}\label{car012}
\begin{array}{rl}
&\ds 2\dbE\int_{Q}\Big(q^{11} v-\lambda\phi\partial _t\th v\Big)dv \\\ns
&\ds
=\dbE\int_{\mathcal{M}}(q^{11}-\lambda\phi\partial _t\th)|v|^2\Big|_{0}^T-\dbE\int_{Q}(\partial_tq^{11}-\lambda\phi\partial _{tt}\th)|v|^2dt\\\ns
&\ds\q
 -\dbE\int_{Q}(q^{11}-\lambda\phi\partial _t\th)|dv|^2\ge -X^1-X^2-X^4.
\end{array}
\end{eqnarray}
%From Lemma 2.6 and Theorem 2.2 it follows that
%$$
%q^{11}=\mathcal{O}_\mu(s^2),\ \partial_t(q^{11})= \th  \mathcal{O}_\mu(s^2).
%$$
%Therefore, by definition of $X^1$, $X^2$, $X^4$ in (\ref{x1}), (\ref{x2}), (\ref{x4}), respectively,  we have
% \begin{eqnarray}\label{car012}
%\begin{array}{rl}
%&\ds2\dbE\int_{Q}q^{11} vdv
%\ge -X^1-X^2-X^4.
%\end{array}
%\end{eqnarray}
Next, by Lemma 2.2 and It\^o's formula, noting that $v=0$ on $\partial Q$,  we have
\begin{eqnarray}\label{car013}
\begin{array}{rl}
&\ds 2\dbE\int_{Q} D_h^2vdv=-2\dbE\int_{Q^*}D_hvD_h(dv)+2\int_{\partial Q} dv t_r(D_hv)\nu\\\ns
&\ds\q\q\q\q\q\q =\dbE\int_{Q^*}|D_h(dv)|^2-\dbE\int_{\mathcal{M}^*}|D_hv|^2\Big|_{0}^T\ge\dbE\int_{Q^*}|D_h(dv)|^2-X^4.
\end{array}
\end{eqnarray}

Further,  by Lemmas 2.1 and 2.2, noting that $v=0$ on $\partial Q$, it yields that
\begin{eqnarray}\label{car014}
\begin{array}{rl}
&\ds h^2\dbE\int_{Q}q^{11} D_h^2v dv\\\ns
&\ds=-h^2\dbE\int_{Q^*}A_hq^{11}D_h(dv)D_hv-h^2\dbE\int_{Q^*}D_hq^{11}A_h(dv)D_hv.
\end{array}
\end{eqnarray}

By $2A_h(dv)D_hv=d(D_h(v^2))-2A_h(dv)D_h(dv)-2A_hvD_h(dv)$,  recalling (\ref{car02}) for the definition of $\Phi(v)$, and  $v=0$ on $\partial Q$ and $t_r(v)=0$ on $\partial Q^*$,  we have
\begin{eqnarray}\label{car015}
\begin{array}{rl}
&\ds-h^2\dbE\int_{Q^*}D_hq^{11}D_hvA_h(dv)+2\dbE\int_{Q}\Phi(v) dv\\\ns
&\ds=-\frac{h^2}{2}\dbE\int_{Q^*}D_hq^{11}d(D_h(v^2))+h^2\dbE\int_{Q^*}D_hq^{11}A_h(dv)D_h(dv)\\
\ns&\ds\q+h^2\dbE\int_{Q^*}D_hq^{11}A_hvD_h(dv)+2\dbE\int_{Q}\Phi(v) dv
\\\ns&\ds=-\frac{h^2}{2}\dbE\int_{M^*}D_hq^{11}D_h(v^2)\Big|_0^T+\frac{h^2}{2}\dbE\int_{Q^*}\partial_t(D _hq^{11})D_h(v^2)dt\\
\ns&\ds\q+\frac{h^2}{2}\dbE\int_{Q^*}D_hq^{11}D_h[(dv)^2].
\end{array}
\end{eqnarray}
Combining (\ref{car014})-(\ref{car015}), by Lemma \ref{l2} and It\^o's formula, noting that $v=0$ on $\partial Q$ and $t_r(v)=0$ on $\partial Q^*$,  we have
\begin{eqnarray}\label{car016}
\begin{array}{rl}
&\ds h^2\dbE\int_{Q}q^{11} D_h^2vdv+2\dbE\int_{Q}\Phi(v)dv\\\ns
&\ds=-\frac{h^2}{2}\[\ds\dbE\int_{\mathcal{M}^*}
A_hq^{11}|D_hv|^2\Big|_0^T-
\dbE\int_{Q^*}\partial_t(A_hq^{11})|D_hv|^2dt -\dbE\int_{Q^*}A_hq^{11}|D_h(dv)|^2\\\ns
&\ds\q-\dbE\int_{\mathcal{M}}D^2_hq^{11}v^2\Big|_0^T+\dbE\int_{Q}\partial_t(D^2_hq^{11})v^2dt
+\dbE\int_{Q}D^2_hq^{11}|dv|^2\].
\end{array}
\end{eqnarray}

 For $\l h (\delta T^2)^{-1}\le 1$,  by (\ref{240107-a}), (\ref{car011})--(\ref{car013}) and  (\ref{car016}), we end up with
\begin{eqnarray}\label{car018}
\begin{array}{rl}
&\ds J_{11}+J_{21}+J_{31}+J_{41}\ge \dbE\int_{Q^*}|D_h(dv)|^2 -X^1-X^2-X^4.
\end{array}
\end{eqnarray}

\ms

\noindent{\bf Step 3. } Let us compute ``$J_{12}+J_{22}+J_{32}+J_{42}$".

Recalling (\ref{car02}) for the definitions of $I_1, I_2$ and $\Phi(v)$, for simplicity, we set
 \bel{24-05-0b}
 q^{12}=r^2 D_h^2\rho A_hD_h\rho,\q q^{22}=r^2 A_h^2\rho A_hD_h\rho, \q q^{32}=\l  \phi\partial_t \th r A_hD_h\rho, \q  q^{42}=rA_hD_h\rho.
 \ee
 By Lemma \ref{l1}, it is easy to see that
  \bel{24-05-0c}
  A_h^2v=v+\frac{h^2}{4}D_h^2v,\q 2D_h^2vA_hD_hv=D_h(|D_hv|^2).
  \ee
 Then, by (\ref{car02}), (\ref{24-05-0b}) and (\ref{24-05-0c}), we know that
  \bel{24-05-0e}\ba{ll}\ds
  J_{12}+J_{22}+J_{32}+J_{42}\\
  \ns\ds=2\ds\dbE\int_Q\[2q^{12}A_h^2v+2q^{22}D_h^2v-2q^{32} v+2q^{42} \Phi(v)\]A_hD_hv dt\\
  \ns\ds=4\ds\dbE\int_Qq^{12}vA_hD_hvdt+\dbE\int_Q\(\frac{h^2}{2}q^{12}+2q^{22}\)D_h(|D_hv|^2)dt\\
  \ns\ds\q-4\ds\dbE\int_Qq^{32}A_hD_hv v dt+4\ds\dbE\int_Qq^{42} A_hD_hv\Phi(v)dt.
  \ea\ee

Noting that $\ds (A_hv)^2=A_h(v^2)-\frac{h^2}{4}(D_hv)^2$, $t_r(v) = 0$ on $\partial{Q^{*}}$ and $v=0$ on $\partial Q$, we have
\begin{eqnarray}\label{car111}
\begin{array}{rl}
&4\ds\dbE\int_Qq^{12}vA_hD_hvdt\\\ns
&=-4\ds\dbE\int_{Q^*}A_hq^{12}D_hvA_hvdt-4\ds\dbE\int_{Q^*} D_hq^{12}|A_hv|^2dt\\\ns
&=-2\ds\dbE\int_{Q^*}A_hq^{12}D_h(v^2)dt-4\dbE\int_{Q^*}  D_hq^{12}A_h(v^2)dt+h^2\dbE\int_{Q^*} D_hq^{12}|D_hv|^2dt\\\ns
&=-2\ds\dbE\int_{Q} D_hA_hq^{12}v^2dt+h^2\dbE\int_{Q^*} D_hq^{12}|D_hv|^2dt.
\end{array}
\end{eqnarray}
Next, by Lemma \ref{l2}, we have
\begin{eqnarray}\label{car112}
\begin{array}{rl}
&\ds\dbE\int_Q\(\frac{h^2}{2}q^{12}+2q^{22}\)D_h(|D_hv|^2)dt\\\ns
&=-\ds\dbE\int_{Q^*}\(\frac{h^2}{2}D_hq^{12}+2D_hq^{22}\)|D_hv|^2dt+\ds\dbE\int_{\partial Q}\(\frac{h^2}{2}q^{12}+2q^{22}\)t_r({|D_hv|^2})\nu dt.
\end{array}
\end{eqnarray}

Combining (\ref{24-05-0b}), (\ref{24-05-0e})--(\ref{car112}) with (\ref{240107-b})--(\ref{240107-c}), by (\ref{x1}) and (\ref{x3}), we have
 \begin{eqnarray}\label{0-car120}
\begin{array}{rl}
&\ds
J_{12}+J_{22}+J_{32}+J_{42}
\ge 6\ds\dbE\int_{Q} s^3\varphi^3\mu^4(\partial_x\psi)^4v^2dt+2\dbE\int_{Q^*}s\varphi \mu^2(\partial_{x}\psi)^2|D_hv|^2dt\\\ns
&\ds\q\q\q\q\q\q\q\q\q\q
 -2\ds\dbE\int_{\partial Q}s\varphi\mu \partial_x \psi  t_r({|D_hv|^2})\nu dt-X^1-X^3\\
\ns &\ds\q\q\q\q\q\q\q\q\q\q-4\ds\dbE\int_Qq^{32}A_hD_hv v dt+4\ds\dbE\int_Qq^{42} A_hD_hv\Phi(v)dt.
\end{array}
\end{eqnarray}
By Lemmas 2.1 and 2.2, and noting that $v=0$ on $\partial Q$, $t_r(v)=0$ on $\partial Q^*$, we have
\begin{eqnarray}\label{1-l6}
\begin{array}{rl}
\ds-4\ds\dbE\int_Q q^{32} A_hD_hv v dt
&= 4\ds\dbE\int_{Q^*} D_hq^{32} |A_hv|^2 dt+2\dbE\int_{Q^*}A_hq^{32}D_h(v^2) dt\\\ns
&\ds = 4\ds\dbE\int_{Q^*} D_hq^{32} |A_hv|^2 dt-2\dbE\int_{Q}D_hA_hq^{32}v^2 dt\\\ns
&\ds \ge -\ds\dbE\int_{Q^*} \th\mathcal{O}_\mu(s^2) |A_hv|^2 dt-2\dbE\int_{Q}\th\mathcal{O}_\mu(s^2)v^2 dt,
\end{array}
\end{eqnarray}
where we have used the following facts:
\begin{eqnarray*}
\begin{array}{rl}
& D_h(q^{32})=D_h(\l  \phi\partial_t \th r A_hD_h\rho)=\th\mathcal{O}_\mu(s^2),\q A_hD_h(q^{32})=\th\mathcal{O}_\mu(s^2).
\end{array}
\end{eqnarray*}
By $|A_hv|^2\le A_h(v^2)$, Lemma 2.2,  and $t_r(v)=0$ on $\partial Q^*$, it yields that
\begin{eqnarray}\label{0-car116}
\begin{array}{rl}
\ds-\dbE\int_{Q^*} \th\mathcal{O}_\mu(s^2) |A_hv|^2 dt&
\ge\ds -\dbE\int_{Q^*} \th\mathcal{O}_\mu(s^2) A_h(v^2) dt= -\dbE\int_{Q^*}\th\mathcal{O}_\mu(s^2) |v|^2 dt.
\end{array}
\end{eqnarray}
Next, by (\ref{car02}), we get
 \begin{eqnarray}\label{car116}
\begin{array}{rl}
&4\ds\dbE\int_Qq^{42} A_hD_hv\Phi(v)dt\\\ns
&\ds=2h^2\ds\dbE\int_Qq^{42} A_hD_hvD_h\[D_h(r D_h^2\rho)A_hv\]dt\\\ns
&\ds=2h^2\dbE\int_Q  q^{42}D^2_h(r D_h^2\rho) A_h^2vA_hD_hv dt+2h^2\dbE\int_Q q^{42} A_hD_h(r D_h^2\rho)|A_hD_hv|^2 dt\\\ns
%\ns&\ds=\dbE\int_Q  s\mathcal{O}_\mu((sh)^2) \(A_h^2vA_hD_hv +|A_hD_hv|^2\) dt\\\ns
&\ds\ge -\dbE\int_Q s\mathcal{O}_\mu((sh)^2)|A_h^2v|^2dt-\dbE\int_Q  s\mathcal{O}_\mu((sh)^2)|A_hD_hv|^2dt,
\end{array}
\end{eqnarray}
where we have used the following facts:
\begin{eqnarray*}
\begin{array}{rl}
& r D^2_h(r D_h^2\rho)A_hD_h\rho=\partial_x ^2 (r\partial_x^2\rho)r\partial_x \rho+s^3\mathcal{O}_\mu((sh)^2)=\mathcal{O}_\mu(s^3),
\end{array}
\end{eqnarray*}
and
\begin{eqnarray*}
\begin{array}{rl}
&r A_hD_h(r D_h^2\rho)A_hD_h\rho=\partial_x (r\partial_x^2\rho)r\partial_x \rho+s^3\mathcal{O}_\mu((sh)^2)=\mathcal{O}_\mu(s^3).
\end{array}
\end{eqnarray*}
By $|A_h^2v|^2\le A_h(|A_hv|^2)$ and Lemma 2.2,  we have
\begin{eqnarray}\label{*}
\begin{array}{rl}
&\ds
 -\dbE\int_{Q}s\mathcal{O}_\mu{ ((sh)^2)}|A^2_hv|^2 dt\ge -\dbE\int_{Q}s\mathcal{O}_\mu{ ((sh)^2)}A_h(|A_hv|^2) dt\\\ns
&\ds= -\dbE\int_{Q^*}s\mathcal{O}_\mu{ ((sh)^2)}|A_hv|^2 dt+\frac{h}{2}\dbE\int_{\partial Q}s\mathcal{O}_\mu{ ((sh)^2)}t_r(|A_hv|^2)dt\\\ns
 &\ds\ge -\dbE\int_{Q^*}s\mathcal{O}_\mu{ ((sh)^2)}|A_hv|^2 dt.
\end{array}
\end{eqnarray}
In the same manner, we can obtain that
\begin{eqnarray}\label{0-car117a}
\begin{array}{rl}
&\ds
  -\dbE\int_{Q^*}s\mathcal{O}_\mu{ ((sh)^2)}|A_hv|^2 dt\ge  -\dbE\int_{Q}s\mathcal{O}_\mu{ ((sh)^2)}|v|^2 dt,
\end{array}
\end{eqnarray}
and
\begin{eqnarray}\label{1-car117a}
\begin{array}{rl}
&\ds
  -\dbE\int_{Q}s\mathcal{O}_\mu{((sh)^2)}|A_hD_hv|^2 dt\ge  -\dbE\int_{Q}s\mathcal{O}_\mu{ ((sh)^2)}|D_hv|^2 dt.
\end{array}
\end{eqnarray}
Then, combining (\ref{1-l6})--(\ref{1-car117a}),  we have
\begin{eqnarray}\label{car117a}
\begin{array}{rl}
&\ds
-4\ds\dbE\int_Qq^{32}A_hD_hv v dt+4\ds\dbE\int_Qq^{42} A_hD_hv\Phi(v)dt\ge -X^1.
\end{array}
\end{eqnarray}
By (\ref{0-car120}) and (\ref{car117a}), we end up with
\begin{eqnarray}\label{car120}
\begin{array}{rl}
&\ds
J_{12}+J_{22}+J_{32}+J_{42}
\ge 6\ds\dbE\int_{Q} s^3\varphi^3\mu^4(\partial_x\psi)^4v^2dt+2\dbE\int_{Q^*}s\varphi \mu^2(\partial_{x}\psi)^2|D_hv|^2dt\\\ns
&\ds\q\q\q\q\q\q\q\q\q\ \ \ 
 -2\ds\dbE\int_{\partial Q}s\varphi\mu \partial_x \psi  t_r({|D_hv|^2})\nu dt-X^1-X^3.
\end{array}
\end{eqnarray}
\medskip

\noindent{\bf Step 4. } Combining (\ref{car018}), (\ref{car120}) with (\ref{car05}), for $\l>1$ and $\lambda h(\delta T^2)^{-1}\le 1$,  it holds that
\begin{eqnarray}\label{car201}
\begin{array}{rl}
 &\ds2\dbE\int_Q I_1I_2\\\ns
 &\ds\ge 6\ds\dbE\int_{Q} s^3\varphi^3\mu^4(\partial_x\psi)^4v^2dt+2\ds\dbE\int_{Q^*}s\varphi\mu^2(\partial_{x}\varphi)^2|D_hv|^2dt+\dbE\int_{Q^*} |D_h(dv)|^2\\\ns
 &\q \ds -2\int_{\partial Q}  s\varphi\mu\partial_{x}\psi t_r(|D_hv|^2)\nu dt-X^1-X^2-X^3-X^4.
 \end{array}
\end{eqnarray}
Next, recall (\ref{car02}) for  the definition of $\Phi(v)$, by  (\ref{240107-a}), we have
\begin{eqnarray}\label{car202}
\begin{array}{rl}
\ds-\dbE\int_Q |\Phi(v)|^2 dt\ge -\dbE\int_Q\mathcal{O}_\mu{ ((sh)^2)}|A_h^2v|^2dt-\dbE\int_Q\mathcal{O}_\mu{ ((sh)^2)}|A_hD_hv|^2dt.
\end{array}
\end{eqnarray}
Similar to the derivation of (\ref{*}), it follows that
\begin{eqnarray}\label{car203}
\begin{array}{rl}
&\ds
-\dbE\int_Q\mathcal{O}_\mu{ ((sh)^2)}|A_h^2v|^2dt\ge -\dbE\int_Q\mathcal{O}_\mu{ ((sh)^2)}|v|^2dt,
\end{array}
\end{eqnarray}
and
\begin{eqnarray}\label{car204}
\begin{array}{rl}
&\ds
-\dbE\int_Q\mathcal{O}_\mu{ ((sh)^2)}|A_hD_hv|^2
\ge -\dbE\int_Q\mathcal{O}_\mu{ ((sh)^2)}|D_hv|^2dt.
\end{array}
\end{eqnarray}
Combining (\ref{car203})-(\ref{car204}) with (\ref{car202}), we have
\begin{eqnarray}\label{car205}
\begin{array}{rl}
&\ds
-\dbE\int_Q |\Phi(v)|^2 dt
\ge -X^1.
\end{array}
\end{eqnarray}
Combining (\ref{car201}), (\ref{car205}) with (\ref{car04}),
\begin{eqnarray*}
\begin{array}{rl}
& \ds 6\ds\dbE\int_{Q} s^3\varphi^3\mu^4(\partial_x\psi)^4v^2dt+2\ds\dbE\int_{Q^*}s\varphi\mu^2(\partial_{x}\varphi)^2|D_hv|^2dt+\dbE\int_{Q^*} |D_h(dv)|^2\\\ns
&\ds\le 2\int_{\partial Q}  s\varphi\mu\partial_{x}\psi t_r(|D_hv|^2)\nu dt+\dbE\int_{Q}|r f|^2dt+X^1+X^2+X^3+X^4.
\end{array}
\end{eqnarray*}

Remember that $\psi$ satisfying (\ref{240107-0a}), then by elementary calculation, we conclude that there exists a positive constant $\mu_0$ such that for any $\mu\ge \mu_0$, one can find two positive constants $\varepsilon_1>0$ and $h_1>0$, with $0<\varepsilon_1\le 1$ and a sufficiently large $\l_1\ge 1$, such that for $\l\ge \l_1$, and $h\le h_1 $ and $\l h (\delta T^2)^{-1}\le \varepsilon_1$, there exist a constant $C>0$ such that
\begin{eqnarray}\label{car11}
\begin{array}{rl}
&\ds \dbE\int_Q s^3\varphi^3\mu^4v^2dt+\ds\dbE\int_{Q^*}s\varphi\mu^2|D_hv|^2dt\\\ns
&\ds\le  \dbE\int_0^T\int_{\omega_0\cap \mathcal{M}}s^3\varphi^3\mu^4v^2dt+\ds\dbE\int_0^T\int_{\omega_0\cap \mathcal{M}^*}s\varphi\mu^2|D_hv|^2dt
\\\ns
&\q\ds+C \(\dbE\int_{Q}|r f|^2dt+\dbE\int_{{Q}}s^2|dv|^2
 \ds+h^{-2}\dbE\int_{\mathcal{M}} |v|^2\Big|_{t=0}
 \ds+h^{-2}\dbE\int_{\mathcal{M}}|v|^2\Big|_{t=T}\),
\end{array}
\end{eqnarray}
where we have used that
$$
|D_hv|^2\le Ch^{-2}(|\tau_{-}v|^2+|\tau_{+}v|^2).
$$
\medskip

\noindent{\bf Step 5. } We now return $v$ to $w$. Recalling $w=\rho v$, we have $D_hw=D_h\rho  A_hv+A_h\rho D_hv$. Therefore, by Lemma 2.3, it holds that
\begin{eqnarray}\label{1-7-01}
\begin{array}{rl}
&\ds\dbE\int_{Q^*}s\varphi\mu^2|r D_hw|^2dt\\\ns
%&\ds \le C\(\dbE\int_{Q^*}s\varphi\mu^2r^2|D_h\rho A_hv|^2dt+\dbE\int_{Q^*}s\varphi\mu^2r^2|A_h\rho D_hv|^2dt\)\\\ns
&\ds\le C\(\dbE\int_{Q^*}s^3\varphi^3\mu^4|A_hv|^2dt+\dbE\int_{Q^*}s^3\mathcal{O}_\mu((sh)^2)|A_hv|^2dt\\\ns
&\ds\qq\q+\dbE\int_{Q^*}s\varphi\mu^2|D_hv|^2dt+\dbE\int_{Q^*}s\mathcal{O}_\mu((sh)^2)|D_hv|^2dt\).
\end{array}
\end{eqnarray}
Similar to the derivation of (\ref{0-car116}), and by Lemma 2.3, we have
\begin{eqnarray}\label{1-7-02a}
\begin{array}{rl}
&\ds\dbE\int_{Q^*}s^3\varphi^3\mu^4|A_hv|^2dt
\le \dbE\int_{Q^*}s^3\varphi^3\mu^4A_h(v^2)dt
= \dbE\int_{Q}s^3\mu^4A_h(\varphi^3)|v|^2dt\\\ns
&\ds= \dbE\int_{Q}s^3\mu^4\varphi^3|v|^2dt+\dbE\int_{Q}s\mathcal{O}_\mu((sh)^2)|v|^2dt.
\end{array}
\end{eqnarray}
And we have
\begin{eqnarray}\label{1-7-02b}
\begin{array}{rl}
&\ds\dbE\int_{Q^*}s^3\mathcal{O}_\mu((sh)^2)|A_hv|^2dt
\le  \dbE\int_{Q}s^3\mathcal{O}_\mu((sh)^2)|v|^2dt.
\end{array}
\end{eqnarray}
 Then, combining (\ref{1-7-01})--(\ref{1-7-02b}), it follows that
\begin{eqnarray}\label{car12}
\begin{array}{rl}
&\ds\dbE\int_{Q^*}s\varphi\mu^2|r D_hw|^2dt\\\ns
&\ds\le C\( \dbE\int_{Q}s^3\mu^4\varphi^3|v|^2dt+\dbE\int_{Q^*}s^3\mathcal{O}_\mu((sh)^2)v^2dt+\dbE\int_{Q^*}s\varphi\mu^2|D_hv|^2dt\\\ns
&\ds\q\q\q+\dbE\int_{Q^*}s\mathcal{O}_\mu((sh)^2)|D_hv|^2dt\).
\end{array}
\end{eqnarray}
Notice that $r A_h \rho D_hv=r D_hw-r D_h\rho A_hv$, then by Lemma 2.3, we have $D_hv=r D_hw-r D_h\rho A_hv-\mathcal{O}_\mu((sh)^2)D_hv$. Then
\begin{eqnarray}\label{1-7-03}
\begin{array}{rl}
&\ds\dbE\int_0^T\int_{\omega_0\cap \mathcal{M}^*}s\varphi\mu^2|D_hv|^2dt\\\ns
&\ds \le C\(\dbE\int_0^T\int_{\omega_0\cap \mathcal{W}^*}s\varphi\mu^2r^2|D_hw|^2dt
+\dbE\int_0^T\int_{\omega_0\cap \mathcal{M}^*} s^3\varphi^3\mu^4 |A_hv|^2dt\\\ns
&\ds\q\q\q+\dbE\int_0^T\int_{\omega_0\cap \mathcal{M}^*} s^3\mathcal{O}_\mu((sh)^2) |A_hv|^2dt
+\dbE\int_0^T\int_{\omega_0\cap \mathcal{M}^*} s\mathcal{O}_\mu((sh)^2) |D_h v|^2dt\).
\end{array}
\end{eqnarray}
Similar to (\ref{*}), and by Lemma 2.3,  we can obtain that
\begin{eqnarray}\label{1-7-04}
\begin{array}{rl}
&\ds\dbE\int_0^T\int_{\omega_0\cap \mathcal{M}^*} s^3\varphi^3\mu^4 |A_hv|^2dt
+\dbE\int_0^T\int_{\omega_0\cap \mathcal{M}^*} s^3\mathcal{O}_\mu((sh)^2) |A_hv|^2dt\\\ns
&\ds\le \dbE\int_{\omega_0\cap \mathcal{M}}s^3\mu^4\varphi^3|v|^2dt+\dbE\int_{\omega_0\cap \mathcal{M}}s^3\mathcal{O}_\mu((sh)^2)|v|^2dt.
\end{array}
\end{eqnarray}
By (\ref{1-7-03})--(\ref{1-7-04}), we have
\begin{eqnarray}\label{car13}
\begin{array}{rl}
&\ds\dbE\int_0^T\int_{\omega_0\cap \mathcal{M}^*}s\varphi\mu^2|D_hv|^2dt\\\ns
&\ds\le C\(\dbE\int_0^T\int_{\omega_0\cap \mathcal{M}^*}s\varphi\mu^2r^2|D_hw|^2dt+\dbE\int_0^T\int_{\omega_0\cap \mathcal{M}} s^3\varphi^3\mu^4 v^2dt\\\ns
&\ds\q\q\q
+\dbE\int_0^T\int_{\omega_0\cap \mathcal{M}} s^3\mathcal{O}_\mu(sh)^2) v^2dt
+\dbE\int_0^T\int_{\omega_0\cap \mathcal{M}^*} s\mathcal{O}_\mu(sh)^2) |D_h v|^2dt\).
\end{array}
\end{eqnarray}
Combining (\ref{car12}), (\ref{car13}) with (\ref{car11}) and noting
$$
\dbE\int_Qs^2|dv|^2=\dbE\int_Qs^2|rdw|^2=\dbE\int_Qs^2|rg|^2dt,
$$
we conclude that, for any $\mu\ge \mu_0$, one can find two positive constants $\varepsilon_2>0$ and $h_2>0$, with $0<\varepsilon_2\le 1$ and a sufficiently large $\l_2\ge 1$, such that for $\l\ge \l_2$, and $h\le h_2 $ and $\l h(\delta T^2)^{-1}\le \varepsilon_2$, there exist a constant $C>0$ such that
\begin{eqnarray}\label{car13a}
\begin{array}{rl}
&\ds\dbE\int_{Q}  s^3\varphi^3\mu^4e^{2s\phi} w^2dt+\ds\dbE\int_{Q^*}s\varphi\mu^2e^{2s\phi} |D_hw|^2dt\\\ns
&\le\ds C\(\dbE\int_0^T\int_{\omega_0\cap \mathcal{W}}s^3\varphi^3\mu^4e^{2s\phi} w^2dt+\ds\dbE\int_0^T\int_{\omega_0\cap \mathcal{W}^*}s\varphi\mu^2e^{2s\phi}|D_hw|^2dt
\\\ns
&\ds\q + \dbE\int_{{Q}}s^2e^{2s\phi}|g|^2dt+\dbE\int_{Q}e^{2s\phi} | f|^2dt
 +h^{-2}\dbE\int_{\mathcal{M}}e^{2s\phi} |w|^2\Big|_{t=0}
+h^{-2}\dbE\int_{\mathcal{M}}e^{2s\phi} |w|^2\Big|_{t=T}\).
\end{array}
\end{eqnarray}
\noindent{\bf Step 6. } Next, we choose a real-valued function $\xi\in C^\infty_0(w;[0,1]) $ such that $\xi = 1$ in $w_0$.
By It\^o's formula and (3.1),  then it is easy to show that
\begin{eqnarray*}
\begin{array}{rl}
& \ds\dbE\int_{\omega\cap \mathcal{M}}s\varphi \xi^2 e^{2s\phi}|w|^2\Big|_{t=T}-\dbE\int_{\omega\cap \mathcal{M}}s\varphi \xi^2  e^{2s\phi}|w|^2\Big|_{t=0}\\\ns
&=\ds\dbE\int_0^T\int_{\omega\cap \mathcal{M}}\xi^2 \[\varphi \partial_t(se^{2s\phi})|w|^2dt+2 s\varphi e^{2s\phi}wdw+ s\varphi  e^{2s\phi}(dw)^2\]\\\ns
&=\ds\dbE\int_0^T\int_{\omega\cap \mathcal{M}}\xi^2 \[\varphi \partial_t(se^{2s\phi})|w|^2+2 s\varphi e^{2s\phi}w(-D_h^2w+f)+ s\varphi  e^{2s\phi}g^2\]dt\\\ns
&=\ds\dbE\int_0^T\int_{\omega\cap \mathcal{M}}\[\xi^2\varphi \partial_t(se^{2s\phi})|w|^2
+2 s\xi^2\varphi e^{2s\phi}wf+ s\xi^2\varphi  e^{2s\phi}g^2\]dt\\\ns
&\ds\q+\ds\dbE\int_0^T\int_{\omega\cap \mathcal{M}^*}\[2 sA_h(\xi^2\varphi e^{2s\phi})|D_hw|^2+2 sD_h(\xi^2 \varphi e^{2s\phi})A_hwD_hw\]dt.
\end{array}
\end{eqnarray*}
Therefore, we conclude that
\begin{eqnarray}\label{car14}
\begin{array}{rl}
&\ds 2\dbE\int_0^T\int_{\omega\cap \mathcal{M}^*} sA_h(\xi^2\varphi e^{2s\phi})|D_hw|^2dt+\dbE\int_0^T\int_{\omega\cap \mathcal{M}}s\xi^2\varphi  e^{2s\phi}g^2dt\\\ns
&\ds=-\dbE\int_0^T\int_{\omega\cap \mathcal{M}}\[\xi^2\varphi \partial_t(se^{2s\phi})|w|^2
+2 s\xi^2\varphi e^{2s\phi}wf\]dt+\dbE\int_{\omega\cap \mathcal{M}}s\varphi \xi^2 e^{2s\phi}|w|^2\Big|_{t=T}\\\ns
&\ds\q-\dbE\int_{\omega\cap \mathcal{M}}s\varphi \xi^2  e^{2s\phi}|w|^2\Big|_{t=0}-2\ds\dbE\int_0^T\int_{\omega\cap \mathcal{W}^*} sD_h(\xi^2 \varphi e^{2s\phi})A_hwD_hwdt.
\end{array}
\end{eqnarray}
By Lemma 2.3, for $\l h(\delta T^2)^{-1}\le 1$, it follows that $A_h(\xi^2\varphi e^{2s\phi})=\xi^2\varphi e^{2s\phi}+\mathcal{O}_\mu((sh)^2)e^{2s\phi}$. Then  the first term in the  the left-hand of (\ref{car14}) can be estimated by
\begin{eqnarray}\label{car15}
\begin{array}{rl}
&\ds 2\dbE\int_0^T\int_{\omega\cap \mathcal{M}^*} sA_h(\xi^2\varphi e^{2s\phi})|D_hw|^2dt\\\ns
&\ds\ge 2\dbE\int_0^T\int_{\omega\cap \mathcal{M}^*} s\xi^2\varphi e^{2s\phi}|D_hw|^2dt-\dbE\int_0^T\int_{\omega\cap \mathcal{M}^*} s\mathcal{O}_\mu((sh)^2) e^{2s\phi}|D_hw|^2dt.
\end{array}
\end{eqnarray}
For the first term in the right-hand of (\ref{car14}), we have
\begin{eqnarray}\label{car16}
\begin{array}{rl}
&\ds -\dbE\int_0^T\int_{\omega\cap\mathcal{M}}\xi^2\varphi \partial_t(se^{2s\phi})|w|^2dt\le \dbE\int_0^T\int_{\omega\cap\mathcal{M}}\th \mathcal{O}_\mu(s^2) e^{2s\phi}|w|^2dt.
\end{array}
\end{eqnarray}
Moreover, it is easy to check that
\begin{eqnarray}\label{car17}
\begin{array}{rl}
 &\ds-\dbE\int_0^T\int_{\omega\cap \mathcal{M}}2 s\xi^2\varphi e^{2s\phi}wfdt\\\ns
 &\ds\le C\dbE\int_0^T\int_{\omega\cap \mathcal{M}}s^2\mu^{2} \varphi^2 e^{2s\phi}|w|^2dt+C\dbE\int_0^T\int_{\mathcal{M}}\mu^{-2} e^{2s\phi}|f|^2dt.
\end{array}
\end{eqnarray}
Next, we estimate  the last term in the  the right-hand of (\ref{car14}). By Lemma 2.3, we have $D_h(\xi^2 \varphi e^{2s\phi})=\partial_x(\xi^2 \varphi e^{2s\phi})+s\mathcal{O}_\mu((sh)^2)e^{2s\phi}=s\mathcal{O}(\mu)\xi\varphi^2e^{2s\phi}+s\mathcal{O}_\mu((sh)^2)e^{2s\phi}$, then
\begin{eqnarray}\label{1-8-01}
\begin{array}{rl}
&\ds -\dbE\int_0^T\int_{\omega\cap \mathcal{M}^*}2 sD_h(\xi^2 \varphi e^{2s\phi})A_hwD_hwdt\\\ns
&\ds =\dbE\int_0^T\int_{\omega\cap \mathcal{M}^*}\Big[s^2\mathcal{O}(\mu)\xi\varphi^2e^{2s\phi}+s\mathcal{O}_\mu((sh)^2)e^{2s\phi}\Big]A_hwD_hwdt\\\ns
&\ds\le \int_0^T\int_{\omega\cap \mathcal{M}^*}\xi^2 s\varphi  e^{2s\phi}|D_hw|^2dt+C\dbE\int_0^T\int_{\omega\cap \mathcal{M}^*} s^3\mu^2\varphi^3  e^{2s\phi}|A_hw|^2dt\\\ns
&\ds\q+\int_0^T\int_{\omega\cap \mathcal{M}^*} s\mathcal{O}_\mu((sh)^2)e^{2s\phi}|A_hw|^2dt
+\int_0^T\int_{\omega\cap \mathcal{M}^*} s\mathcal{O}_\mu((sh)^2)e^{2s\phi}|D_hw|^2dt.
\end{array}
\end{eqnarray}
Similar to (\ref{*}), and noting that $A_h(\varphi^3  e^{2s\phi})=\varphi^3  e^{2s\phi}+ \mathcal{O}_\mu((sh)^2) e^{2s\phi}$,  $A_h( e^{2s\phi})= \mathcal{O}_\mu(1)e^{2s\phi}$, we have
\begin{eqnarray}\label{1-8-02}
\begin{array}{rl}
&\ds \dbE\int_0^T\int_{\omega\cap \mathcal{M}^*} s^3\mu^2\varphi^3  e^{2s\phi}|A_hw|^2dt+\int_0^T\int_{\omega\cap \mathcal{M}^*}s\mathcal{O}_\mu((sh)^2)e^{2s\phi}|A_hw|^2dt\\\ns
%&\ds \le\dbE\int_0^T\int_{\omega\cap \mathcal{M}^*} s^3\mu^2\varphi^3  e^{2s\phi}A_h(w^2)dt\\\ns
&\ds\le\dbE\int_0^T\int_{\omega\cap \mathcal{M}} s^3\mu^2\varphi^3  e^{2s\phi}|w|^2dt
+\int_0^T\int_{\omega\cap \mathcal{M}}s^3\mathcal{O}_\mu((sh)^2)e^{2s\phi}|w|^2dt.
\end{array}
\end{eqnarray}
Hence, combining (\ref{1-8-01}) and (\ref{1-8-02}), the last term in the  the left-hand of (\ref{car14}) can be estimated by
\begin{eqnarray}\label{car18}
\begin{array}{rl}
&\ds -\dbE\int_0^T\int_{\omega\cap \mathcal{M}^*}2 sD_h(\xi^2 \varphi e^{2s\phi})A_hwD_hwdt\\\ns
&\ds\le \int_0^T\int_{\omega\cap \mathcal{M}^*}\xi^2 s\varphi  e^{2s\phi}|D_hw|^2dt+C\dbE\int_0^T\int_{\omega\cap \mathcal{M}^*} s^3\mu^2\varphi^3  e^{2s\phi}|w|^2dt\\\ns
&\ds\q+\int_0^T\int_{\omega\cap \mathcal{M}^*} s^3\mathcal{O}_\mu((sh)^2)e^{2s\phi}|w|^2dt
+\int_0^T\int_{\omega\cap \mathcal{M}^*} s\mathcal{O}_\mu((sh)^2)e^{2s\phi}|D_hw|^2dt.
\end{array}
\end{eqnarray}
Therefore, combining (\ref{car14})-(\ref{car17}), (\ref{car18}), we conclude that for any $\mu\ge \mu_0$, one can find two positive constants $\varepsilon_3>0$ and $h_3>0$, with $0<\varepsilon_3\le 1$ and a sufficiently large $\l_3\ge 1$, such that for $\l\ge \l_3$, and $h\le h_3 $ and $\l h (\delta T^2)^{-1}\le \varepsilon_3$, it holds that
\begin{eqnarray}\label{car19}
\begin{array}{rl}
&\ds \dbE\int_0^T\int_{\omega_0\cap \mathcal{M}^*} s\varphi e^{2s\phi}|D_hw|^2dt
\le \dbE\int_0^T\int_{\omega\cap \mathcal{M}^*} s\xi^2\varphi e^{2s\phi}|D_hw|^2dt\\\ns
&\ds\le C\(\dbE\int_0^T\int_{\omega\cap \mathcal{M}} s^3\mu^2\varphi^3  e^{2s\phi}|w|^2dt+\dbE\int_0^T\int_{\omega\cap \mathcal{M}}\mu^{-2} e^{2s\phi}|f|^2dt\\\ns
&\ds\q \q\q+h^{-1}\dbE\int_{\omega\cap \mathcal{M}} e^{2s\phi}|w|^2\Big|_{t=T}\).
\end{array}
\end{eqnarray}
By (\ref{car13a}) and  (\ref{car19}),  then
we complete the proof.
\endpf

\ms

\section{The observability for the  backward stochastic semi-discrete  parabolic equations}
 The controllability result for (\ref{e*03}) can be established by proving an observability estimate for its adjoint system. In this section, we first consider the following backward stochastic semi-discrete  parabolic equation:
\begin{eqnarray}\label{e*03b1}\left\{
\begin{array}{lll}
\ds \mbox{d}z+ D_h^2zdt
=-(a_1 z+a_zZ)dt+ZdB(t),
\  (x,t)\in Q,\\[3mm]
z(0,t)=0, \ z(1,t)=0,\ \  t\in (0,T),\\[3mm]
 z(x,T)=z_T, \ x\in \mathcal{M}.
\end{array}
\right.
\end{eqnarray}
Then the following observability inequality holds.
\bt
There exist  $h_0>0$ and  $C$ such that for all $h\le h_0$,   the following observability inequality holds
\begin{eqnarray*}\
\begin{array}{rl}
&\mathbb{E}\displaystyle\int_{\mathcal{M}}  |z(0)|^2
 \le
C\(\dbE\int_{Q}| Z|^2dt+\dbE\int_{0}^T\int_{\omega\cap\mathcal{M}} |z|^2dt
+ e^{-\frac{C}{h}}\dbE\int_{\mathcal{M}} |z|^2\Big|_{t=T}\).
 \end{array}
\end{eqnarray*}
for any solution to system  $(\ref{e*03b1})$.
\et
{\bf Proof.} Let
$$
A=|a_1|_{L_{\dbF}^{\infty}(0, T; L_h^{\infty}(\mathcal{M}))}+|a_2|_{L_{\dbF}^{\infty}(0, T; L_h^{\infty}(\mathcal{M}))}.
$$
By Theorem 3.1, it deduce that
\begin{eqnarray*}\
\begin{array}{rl}
&\ds\dbE\int_{Q} s^3e^{2s\phi} z^2dt+\ds\dbE\int_{Q^*}se^{2s\phi} |D_hz|^2dt\\\ns
&\ds\le C\[\dbE\int_{\omega\cap\mathcal{M}} s^3e^{2s\phi} z^2dt+\dbE\int_{Q}s^2 e^{2s\psi} | Z|^2dt+A^2\dbE\int_{Q} e^{2s\phi}(|z|^2+ | Z|^2)dt
 \ds\\\ns
 &\q\q\q\ds+h^{-2}\dbE\int_{\mathcal{M}}e^{2s\phi} |z|^2\Big|_{t=0}
+h^{-2}\dbE\int_{\mathcal{M}}e^{2s\phi} |z|^2\Big|_{t=T}\].
\end{array}
\end{eqnarray*}
After choosing a sufficiently large $\l$, we fix parameters $\l$. It follows that
\begin{eqnarray*}\
\begin{array}{rl}
&\ds\dbE\int_{Q} \th^3 e^{2s\phi} z^2dt+\ds\dbE\int_{Q^*}\th e^{2s\phi} |D_hz|^2dt\\\ns
&\ds\le C\[\dbE\int_{\omega\cap\mathcal{M}} \th^3e^{2s\phi} z^2dt+\dbE\int_{Q}\th^2e^{2s\phi} | Z|^2dt
 +h^{-2}\dbE\int_{\mathcal{M}}e^{2s\phi} |z|^2\Big|_{t=0}\\\ns
  &\q\q\q\ds
 \ds+h^{-2}\dbE\int_{\mathcal{M}}e^{-2s\phi} |z|^2\Big|_{t=T}
 \].
\end{array}
\end{eqnarray*}
Recalling the definition of $\th$, we have $\th(t)\ge \th(\frac{T}{2})\ge T^{-2}$, since $0<\delta<\frac{1}{2}$, and $\th\le \th(\frac{T}{4})\le \frac{16}{3T^2}$  for $t\in [\frac{T}{4},\frac{3T}{4}]$. Then
\begin{eqnarray}\label{2-13-01}
\begin{array}{rl}
&\ds\dbE\int_{Q} \th^3 e^{2s\phi} z^2dt\ge \dbE\int_{\frac{T}{4}}^{\frac{3T}{4}}\int_{\mathcal{M}} \th^3 e^{2s\phi} z^2dt
 \ge Ce^{-CT^{-2}}\dbE\int_{\frac{T}{4}}^{\frac{3T}{4}}\int_{\mathcal{M}}  z^2dt.
\end{array}
\end{eqnarray}
On the other hand, notice that $d(|z|^2)=2zd{z}+|dz|^2$.  Hence,  for any $0\leq t_1\leq t_2\leq T$,
\begin{eqnarray*}
&&\mathbb{E}\int_{\mathcal{M}} |z(t_2)|^2-\mathbb{E}\int_{\mathcal{M}} |z(t_1)|^2\\
&&=\mathbb{E}\int^{t_2}_{t_1}\int_ {\mathcal{M}} \Big\{
2z\Big[ -D_h^2z
-(a_1 z+a_zZ)\Big]+|Z|^2 \Big\}dt\\
&&=\mathbb{E}\int_ {t_1}^{t_2}\int_{\mathcal{M}}  \Big[
2|D_h z|^2-2a_1|z|^2-2a_2zZ+|Z|^2\Big]dt\\
&&\ge -\mathbb{E}\int_{t_1}^{t_2}\int_ {\mathcal{M}}
\Big[C\sqrt{A}|z|^2+CA|z|^2\Big]dxdt\\
&&\ge -C(1+A)\mathbb{E}\int_{t_1}^{t_2}\int_ {\mathcal{M}} |z|^2dt.
\end{eqnarray*}
By the Gronwall inequality, it holds that
\begin{eqnarray}\label{2-13-02}
\mathbb{E}\displaystyle\int_{\mathcal{M}}  |z(t_1)|^2
\leq C\mathbb{E}\displaystyle\int_{\mathcal{M}}  |z(t_2)|^2, \quad \ 0\leq t_1\leq t_2\leq T.
\end{eqnarray}
 Noting that  $\th(T)=\th(0)\ge \frac{2}{3}({\delta T^2})^{-1}$, then
\begin{eqnarray}\label{2-13-03}
\begin{array}{rl}
&\ds h^{-2}\dbE\int_{\mathcal{M}}e^{2s\phi} |z|^2\Big|_{t=T}+h^{-2}\dbE\int_{\mathcal{M}}e^{2s\phi} |z|^2\Big|_{t=T}\\\ns
&\ds\le C h^{-2}e^{-\frac{C}{\delta T^2}}\[\dbE\int_{\mathcal{M}} |z|^2\Big|_{t=0}
 +\dbE\int_{\mathcal{M}} |z|^2\Big|_{t=T}\]\\\ns
 &\ds\le C h^{-2}e^{-\frac{C}{\delta T^2}}\dbE\int_{\mathcal{M}} |z|^2\Big|_{t=T}.
\end{array}
\end{eqnarray}
Moreover, $\th(t)\ge \th(\frac{T}{2})\ge{ T^{-2}}$, we conclude that
\begin{eqnarray}\label{2-13-04}
\begin{array}{rl}
&\ds\dbE\int_{\omega\cap\mathcal{M}} \th^3e^{2s\phi} z^2dt+\dbE\int_{Q}\th^2e^{2s\phi} | Z|^2dt\\\ns
 &\ds\le  e^{-\frac{C}{ T^2}}\dbE\int_{\omega\cap\mathcal{M}}z^2dt+e^{-\frac{C}{ T^2}}\dbE\int_{Q} | Z|^2dt.
\end{array}
\end{eqnarray}
Combining  (\ref{2-13-01})--(\ref{2-13-04}), we have 
\begin{eqnarray*}\
\begin{array}{rl}
&\mathbb{E}\displaystyle\int_{\mathcal{M}}  |z|^2(0)\leq CT\mathbb{E}\int^{\frac{3T}{4}}_{\frac{T}{4}}
\int_{\mathcal{M}} |z|^2(x, t)dt\leq  C  e^{CT^{-2}}\dbE\int_{Q} \th^3 e^{-2s\varphi} z^2dt\\\ns
 &\ds\le
C\(\dbE\int_{\omega\cap\mathcal{M}}z^2dt+\dbE\int_{Q} | Z|^2dt\)
 + C h^{-2}e^{-(\frac{C}{\delta T^2}-\frac{C}{T^2})}\dbE\int_{\mathcal{M}} |z|^2\Big|_{t=T}\\\ns
  &\ds\le C\(\dbE\int_{\omega\cap\mathcal{M}}z^2dt+\dbE\int_{Q} | Z|^2dt\)
 + C h^{-2}e^{-\frac{C}{\delta T^2}}\dbE\int_{\mathcal{M}} |z|^2\Big|_{t=T},
 \end{array}
\end{eqnarray*}
 for $\delta\le \delta_0$, where  $\delta_0$ is sufficiently small.
Notice that we need $h<h_0$ and $\l h(\delta T^2)^{-1}\le \eps_0,$ where $\l$ is fixed, then let
$h_1= \frac{\eps_0 \delta_0 T^2}{\l}$. For every  $h\le \min \{h_0, h_1\}$, we choose $\delta=\frac{h}{h_1}\delta_0\le\delta_0$, then  $\l h(\delta T^2)^{-1}=\eps_0$. Therefore,
\begin{eqnarray*}\
\begin{array}{rl}
&\mathbb{E}\displaystyle\int_{\mathcal{M}}  |z|^2(0)
 \le
C\(\dbE\int_{\omega\cap\mathcal{M}}z^2dt+\dbE\int_{Q} | Z|^2dt\)
 + C h^{-2}e^{-\frac{C}{h}}\dbE\int_{\mathcal{M}} |z|^2\Big|_{t=T}\).
 \end{array}
\end{eqnarray*}
This
completes the proof of Theorem 4.1.
\endpf

\section{The controllability for the stochastic semi-discrete  parabolic equation}
This section is devoted to the proof of our
controllability result, Theorem 1.1. The proof is almost standard. However,
we give the details for the sake of completeness.

\noindent{\bf  Proof of Theorem 1.1.}
 For any $y_0\in L^2_h(\mathcal{M})$, define a  functional $J(\cdot)$ on $L^2_{\mathcal{F}_T}(\Omega,L^2_h(\mathcal{M}))$ as follows:
$$
J(z_T)=\frac{1}{2}\dbE\int_{Q}| Z|^2dt+\frac{1}{2}\dbE\int_{0}^T\int_{\omega\cap\mathcal{M}}|z|^2 dt+\frac{\eps}{2}\dbE\Vert z_T\Vert^2_{L^2_h(\mathcal{M})}-\dbE\langle y_0, z(0)\rangle_{L^2_h(\mathcal{M})},
$$
where $\eps= e^{-\frac{C}{h}}$, and $(z, Z)$ is the solution to (\ref{e*03b1}) with final value $z_T$.
 It is a simple matter to see that the functional $J(\cdot)$ is continuous and convex.  Now we prove that it is coercive. By Young's inequality and Theorem 4.1, we have
\begin{eqnarray*}
\begin{array}{rl}
&\ds J(z_T)\ge\frac{1}{2}\dbE\int_{Q}| Z|^2dt+\frac{1}{2}\dbE\int_{0}^T\int_{\omega\cap\mathcal{M}}|z|^2 dt+\frac{\eps}{2}\Vert z_T\Vert^2_{L^2_h(\mathcal{M})}\\\ns
&\ds\q\q \q\q-\varepsilon\Vert z(0)\Vert^2_{L^2_h(\mathcal{M})}-\frac{1}{\varepsilon}\Vert y_0\Vert^2_{L^2_h(\mathcal{M})}\\\ns
&\ds\q\q\ \ge \frac{1}{4}\dbE\int_{Q}| Z|^2dt+\frac{1}{4}\dbE\int_{0}^T\int_{\omega\cap\mathcal{M}}|z|^2 dt+\frac{\eps}{4}\Vert z_T\Vert^2_{L^2_h(\mathcal{M})}-{C}\Vert y_0\Vert^2_{L^2_h(\mathcal{M})}.
 \end{array}
\end{eqnarray*}
Then, $J(\cdot)$ is coercive.  Therefore,  $J(\cdot)$ admits a unique minimizer $z_T^*$. Denote $(z^*, Z^*)$ the solution to (\ref{e*03b1}) with final value $z_T^*$. Then we have
\begin{eqnarray*}
\begin{array}{rl}
\ds \dbE\int_{Q}ZZ^*dt+\dbE\int_{0}^T\int_{\omega\cap\mathcal{M}}zz^* dt+\eps \dbE \langle z_T,  z_T^*\rangle_{L^2_h(\mathcal{M})}-\dbE\langle y_0, z(0)\rangle_{L^2_h(\mathcal{M})}=0,
 \end{array}
\end{eqnarray*}
for any $z_T$, with the associated solution $(z, Z)$ to (\ref{e*03b1}). We set $u^*= -z^*\chi_{\omega\cap\mathcal{M}}$ and $v^*=-Z^*$. We claim that (\ref{e*03b1}) with $u=u^*$ and $v=v^*$ satisfy our needs.
By the duality of  (\ref{e*03b1}) with  (\ref{e*03a}), we have
\begin{eqnarray*}
\begin{array}{rl}
\ds \dbE \int_{\mathcal{M}}y(T)z_T- \dbE\int_{\mathcal{M}}y_0z(0)=\dbE\int_{Q}vZdt+\dbE\int_{0}^T\int_{\omega\cap\mathcal{M}}zu dt.
 \end{array}
\end{eqnarray*}
Then $y(T)=\eps z_T^*$. On the other hand,
\begin{eqnarray*}
\begin{array}{rl}
&\ds \dbE\int_{Q}|Z^*|^2dt+\dbE\int_{0}^T\int_{\omega\cap\mathcal{M}}|z^*|^2 dt+\eps \dbE \Vert z_T^*\Vert^2_{L^2_h(\mathcal{M})}= \dbE\langle y_0, z^*(0)\rangle_{L^2_h(\mathcal{M})}\\\ns
&\ds\ds \le  \frac{1}{4\varepsilon}\dbE\int_{\mathcal{M}}|y_0|^2dt+\varepsilon\dbE\int_{\mathcal{M}}|z^*(0)|^2dt.
 \end{array}
\end{eqnarray*}
Notice that
\begin{eqnarray*}
\begin{array}{rl}
&\ds\dbE\int_{\mathcal{M}}|z^*(0)|^2dt\le{C}\(\dbE\int_{Q}| Z^*|^2dt+\dbE\int_{0}^T\int_{\omega\cap\mathcal{M}}|z^*|^2 dt
 + \eps \dbE\int_{\mathcal{M}} |z^*|^2\Bigg|_{t=T}\).
 \end{array}
\end{eqnarray*}
Therefore
 \begin{eqnarray*}
\begin{array}{rl}
&\ds \dbE\int_{Q}|Z^*|^2dt+\dbE\int_{0}^T\int_{\omega\cap\mathcal{M}}|z^*|^2 dt
 \le  C\dbE\int_{\mathcal{M}}|y_0|^2dt,
 \end{array}
\end{eqnarray*}
and
 \begin{eqnarray*}
\begin{array}{rl}
&\ds\eps \dbE \Vert z_T^*\Vert^2_{L^2_h(\mathcal{M})}
 \le  C\dbE\int_{\mathcal{M}}|y_0|^2dt.
 \end{array}
\end{eqnarray*}
Hence
 \begin{eqnarray*}
\begin{array}{rl}
&\ds \dbE\int_{Q}|v|^2dt+\dbE\int_{0}^T\int_{\omega\cap\mathcal{M}}|u|^2 dt
 \le  C\dbE\int_{\mathcal{M}}|y_0|^2dt,
 \end{array}
\end{eqnarray*}
and
\begin{eqnarray*}
\begin{array}{rl}
&\ds \dbE \Vert y(T)\Vert^2_{L^2_h(\mathcal{M})}
 \le  C\eps\dbE\int_{\mathcal{M}}|y_0|^2dt.
 \end{array}
\end{eqnarray*}
Then, we complete the proof.
\endpf
%%%%%%%%%%%%%

\end{document}